\documentclass[11pt]{amsart}
\usepackage{amsfonts,amssymb, graphicx,a4}
\usepackage{color}\definecolor{Red}{cmyk}{0,1,1,0}
\usepackage{mathrsfs}
\usepackage{textcomp}
\usepackage{enumitem}
\usepackage{tikz, tkz-euclide}
\usetikzlibrary{calc}
\usetikzlibrary{decorations.markings}
\usepackage{ulem}

\DeclareMathOperator*{\argmin}{arg\,min}

\numberwithin{equation}{section}

\usepackage{bbm}
  \newtheorem*{theorem*}        {Theorem}
	\newtheorem*{conjecture*}   {Conjecture}
  \newtheorem{theorem}           {Theorem}
  \newtheorem{lemma}              {Lemma}
  \newtheorem*{lemma*}          {Lemma}
    \newtheorem*{claim*}          {Claim}
  \newtheorem{definition}         {Definition}
  \newtheorem{corollary}          {Corollary}
  \newtheorem{proposition}      {Proposition}

\begin{document}


\voffset=-1.5truecm\hsize=16.5truecm    \vsize=24.truecm
\baselineskip=14pt plus0.1pt minus0.1pt \parindent=12pt
\lineskip=4pt\lineskiplimit=0.1pt      \parskip=0.1pt plus1pt

\def\ds{\displaystyle}\def\st{\scriptstyle}\def\sst{\scriptscriptstyle}

\global\newcount\numsec\global\newcount\numfor
\gdef\profonditastruttura{\dp\strutbox}
\def\senondefinito#1{\expandafter\ifx\csname#1\endcsname\relax}
\def\SIA #1,#2,#3 {\senondefinito{#1#2}
\expandafter\xdef\csname #1#2\endcsname{#3} \else
\write16{???? il simbolo #2 e' gia' stato definito !!!!} \fi}
\def\etichetta(#1){(\veroparagrafo.\veraformula)
\SIA e,#1,(\veroparagrafo.\veraformula)
 \global\advance\numfor by 1
 \write16{ EQ \equ(#1) ha simbolo #1 }}
\def\etichettaa(#1){(A\veroparagrafo.\veraformula)
 \SIA e,#1,(A\veroparagrafo.\veraformula)
 \global\advance\numfor by 1\write16{ EQ \equ(#1) ha simbolo #1 }}
\def\BOZZA{\def\alato(##1){
 {\vtop to \profonditastruttura{\baselineskip
 \profonditastruttura\vss
 \rlap{\kern-\hsize\kern-1.2truecm{$\scriptstyle##1$}}}}}}
\def\alato(#1){}
\def\veroparagrafo{\number\numsec}\def\veraformula{\number\numfor}
\def\Eq(#1){\eqno{\etichetta(#1)\alato(#1)}}
\def\eq(#1){\etichetta(#1)\alato(#1)}
\def\Eqa(#1){\eqno{\etichettaa(#1)\alato(#1)}}
\def\eqa(#1){\etichettaa(#1)\alato(#1)}
\def\equ(#1){\senondefinito{e#1}$\clubsuit$#1\else\csname e#1\endcsname\fi}
\let\EQ=\Eq


\def\\{\noindent}
\def\v{\vskip.1cm}
\def\vv{\vskip.2cm}


%
  \def\P{\mathop{\textrm{\rm P}}\nolimits}                  
  \def\d{\mathop{\textrm{\rm d}}\nolimits}                  
  \def\exp{\mathop{\textrm{\rm exp}}\nolimits}              
	\def\supp{\mathop{\textrm{\rm supp}}\nolimits}            
	\def\Int{\mathop{\textrm{\rm Int}}\nolimits}            
	\def\Ext{\mathop{\textrm{\rm Ext}}\nolimits}            
	\def\LRO{\mathop{\textrm{\rm LRO}}\nolimits}            
    \def\sf{\mathop{\textrm{\rm sf}}\nolimits}            
    \def\conv{\mathop{\textrm{\rm conv}}\nolimits}            

    \newcommand\bfblue[1]{\textcolor{blue}{\textbf{#1}}}
\newcommand\blue[1]{\textcolor{blue}{}}

\thispagestyle{empty}

\begin{center}
{\LARGE Entropy Compression Method and Legitimate Colorings in Projective Planes}
\vskip.5cm
Rodrigo Bissacot and Lu\'{i}s Doin
\vskip.3cm
\begin{footnotesize}
Institute of Mathematics and Statistics (IME-USP), University of S\~{a}o Paulo, Brazil\\

\end{footnotesize}
\vskip.1cm
\begin{scriptsize}
emails: rodrigo.bissacot@gmail.com; luisdoin@gmail.com
\end{scriptsize}

\end{center}

\def\be{\begin{equation}}
\def\ee{\end{equation}}

\vskip1.0cm
\begin{quote}
{\small

\textbf{Abstract.} \begin{footnotesize} We prove that the entropy compression method systematized by L. Esperet and A. Parreau \cite{Entropy Compression} can be applied to any problem formulated in the variable version of the Lov\'asz Local Lemma. As an application, we prove the existence of legitimate colorings for projective planes with small orders, which extends results of N. Alon and Z. F\"uredi \cite{Alon}. In fact, we allow different numbers of colors, proving that projective planes of any order can be legitimate colored with 42 colors.  \end{footnotesize}

}
\end{quote}


{\noindent \footnotesize{\bf Keywords:} Finite Projective Plane, Entropy Compression, Lov\'asz Local Lemma, Probabilistic Method.}

\vskip.1cm

{\noindent \footnotesize {\bf Mathematics Subject Classification (2000):} 05-XX, 05C15, 05D40.}

\section{Introduction}

As is well known, the Lov\'asz Local Lemma (LLL, for short) is a powerful tool within the Probabilistic Method for existence proofs. The downside being that it doesn't build the object for which it proves the existence. The celebrated result obtained by Moser and Tardos \cite{MT} in 2010 solved this issue. For a special case of the LLL, called \textit{variable version} - this name is due to Kolipaka and Szegedy \cite{MT2} - the Moser-Tardos algorithm searches in the probability space for a point that avoids all bad events. To prove that their algorithm halts, Moser and Tardos created a new type of argument based on information theory, more precisely on entropy compression, called \textit{entropy compression method} - this name is due to Tao \cite{Tao}. In 2013, Grytczuk, Kozik and Micek \cite{GKeM} obtained new results for nonrepetitive sequences using an entropy compression argument. In the same year Dujmovi\'{c}, Joret, Kozik and Wood \cite{Entropy1} extended their arguments for graphs of bounded degree. Also in 2013,  Esperet and Parreau \cite{Entropy Compression} systematized this approach for a certain group of problems. Many results were recently obtained applying this last method as in \cite{Entropy3, Entropy4, Entropy Compression, Entropy2, Entropy5, Entropy6}.

The entropy compression method schematized in \cite{Entropy Compression} aims to find a coloring $C$ of a graph $G$ such that for any forbidden \textit{configuration} $(H,c)$ and any copy $H'$ of $H$ in $G$, the restriction of the coloring $C$ to $H'$ is not congruent to $c$. Here, by a configuration we mean a pair $(H, c)$ of a subgraph $H$ of $G$ and its coloring $c$. We say that two colorings of the same graph are congruent if one can be obtained from the other by a permutation of the color names. They showed that the approach can be applied to any vertex-coloring (or edge-coloring) problem in a graph $G$ where:
\vskip.1cm
(i) The problem can be stated as to find a coloring where some configurations are forbidden.
\vskip.1cm
(ii) For any configuration $(H,c)$ and any vertex $v$ of $H$, there exist $k_{(H,c)}$ fixed vertices different from $v$ in $H$ for which, if we know their colors, there is a unique way to extend this partial coloring to a coloring of $H$ congruent to $c$. 
\vskip.1cm
   In the present paper we eliminate the second constraint showing that the entropy compression method can be applied to any problem where we are seeking a coloring where some configurations are forbidden. As a consequence, we have the two main contributions of this note: \\
  \vskip.1cm
(1) This slight modification of the Esperet-Parreau method allows us to apply their approach to problems stated in the variable version of the LLL.\\
\vskip.1cm
(2) Applying this less restrictive entropy compression version, it is proved the existence of a legitimate $8$-coloring for projective planes of order greater than $10^{54}$. The previous result was obtained for planes of order at least $10^{250}$ in \cite{Alon}. Further, we have results for different number of colors which indicate the existence of legimate $8$-colorings for any projective plane and not only for planes of large orders. In particular, with 42 colors any projective plane has a legitimate coloring.

Some other attempts to prove a natural connection between LLL type results and entropy compression have been made in \cite{Aldo} and \cite{Anton}, it seems that this connection was clarified recently in \cite{AIS}. In fact, Bernshteyn \cite{Anton} proved a theorem which recovers the same results as the entropy compression in some examples. Besides, his theorem implies the standard lopsided Lov\'{a}sz Local Lemma but there is no formal proof that the result implies the entropy compression results in general and Esperet-Parreau theorem \cite{Entropy Compression}. Also, there are examples where improvements of LLL (via statistical mechanics) give the same answer as the entropy compression, for instance, see the example of the independent set problem in \cite{Aldo, BFPS}. Even the original lopsided LLL was able to achieve the same result as the entropy compression in \cite{Anton1} and surpass in \cite{GST}.

This text is organized as follows: in section 2 we prove that the entropy compression method can be applied to any problem in the variable version of the Lov\'asz Local Lemma. As a corollary, in section 3 we study the problem of legitimate coloring of finite projective planes introduced by Alon and F\"uredi in \cite{Alon}. 

 \\

\textbf{Warning! } We assume that the reader is familiar with the variable version of the Lov\'asz Local Lemma. In adition, we do not repete proofs of the main references \cite{Entropy Compression} and \cite{Alon}. Our main goal is to show that a small modification in the algorithm allows us to use this strategy in all examples formulated in the variable version that welcomes an uniform treatment.
  
\section{Variable Version via Entropy Compression}  

In this section we start with some definitions and fix our setting.

\medskip 

In this paper all the graphs are finite, simple and undirected. Let $G$ be a graph (or a hypergraph) with vertex set $V$ for which we are seeking a certain vertex-coloring such that there are \textit{bad settings} to avoid. A bad setting is a pair $(Z, c)$ where $Z \subset V$ and $c$ is a given coloring of $Z$. Let:

\vskip.1cm

1. $\mathcal{Z} := \{ Z \subset V : \text{There exists} \ c \text{ such that } (Z,c) \text{ is} $ $\text{a bad setting} \} \text{ indexed by a finite}\\ \text{\ set }  I \subset \mathbb{N}$.

\vskip.1cm

2. For each $i \in I$: 

\begin{minipage}[t]{0.063\textwidth}
$\quad \ \ $a.
\end{minipage}\begin{minipage}[t]{\textwidth}
$C_i := \{ c :  (Z_i, c) \ \text{is a bad}$ $\text{setting} \}$. \end{minipage}
\vskip.1cm

$\quad \ \ $b. $\mathscr{Z}_i := \{ (Z_i, c) : c \in C_i\}$.

\vskip.1cm

\begin{minipage}[t]{0.063\textwidth}
$\quad \ \ $c. \end{minipage}\begin{minipage}[t]{0.83\textwidth}For any vertex $v$ of $Z_i$, take $k_i$ fixed vertices different from $v$  in $Z_i$ for which, if we know their colors, there are at most $m_i$ ways to extend this partial coloring to a coloring $c \in C_i$. This is always possible since we will always work with a finite number of colors.\end{minipage}

\vskip.1cm

$\quad \ \ $d. $l_i := \vert Z_i\vert - k_i.$

\vskip.1cm

3. $E := \{ l \in \mathbb{N} : \exists \; i \in I \text{ such that } l_i = l\}$.

\vskip.1cm
\begin{minipage}[t]{0.027\textwidth}	
4. \end{minipage}\begin{minipage}[t]{0.866\textwidth}For each $l \in E$ and $v \in V$, let $Q_l(v)$ be the set of subsets of $V$ containing $v$ and equal to some $Z_j \in \mathcal{Z}$ with $l_j = l$. Then $d_l := \max_{v \in V} \vert Q_l(v) \vert.$
\end{minipage}

\subsection{The Algorithm}

\

\

Now we introduce the setting of the variable version of the LLL and apply the entropy compression algorithm on it.  For precise definitions see the seminal paper \cite{MT}. Let $\Psi = \{\psi_1, \ldots, \psi_n\}$ be a set of independent discrete random variables taking values on a finite set and $\mathcal{A} = \{ A_x : x \in X \}$  be a family of (bad) events indexed by a finite set $X$ and  determined by such variables. To each $x \in X$, let $c(A_x)$ be the set of configurations of the variables of vbl($A_x$)  that makes $A_x$ occur. In order to apply the entropy compression method we need to define the bad settings. Let $G$ be a graph with $V(G) = \Psi$. The problem of finding a configuration of $\Psi$ which avoids all events in $\mathcal{A}$ can be seen as a problem of finding a coloring of the vertices of $G$ in the sense that, given $i \in [n] := \{1, 2, \dots, n\}$, when  $\psi_i$ is sampled, the value sampled is its color. Then 
$\mathcal{Z} = \{ \text{vbl}(A_x) : A_x \in \mathcal{A}\}, \ \ I = X \text{ and, for each } x \in X, \  \mathscr{Z}_x = \{ (\text{vbl}(A_x),c) : c \in c(A_x)\}$. Now, let $\Psi_1, \ldots, \Psi_n$  be vectors that store $t$ samples (for some large $t$) of the random variables $\psi_1, \ldots, \psi_n$, respectively, and $x_{\psi_1},\ldots, x_{\psi_n}$ variables that can assume values in $Im(\psi_1) \cup \{ \circ \}, \ldots, Im(\psi_n) \cup \{ \circ  \} $, respectively. The symbol $\circ$ indicates the situation when the variable is uncolored and we define this as the initial value of the variables $x_{\psi_1},\ldots, x_{\psi_n}$. For each~ $j \in [n]$ and each~ $l \in E$, fix an order on the set $\{ \text{vbl}(A_x) : x \in X, \ \psi_j \in \text{vbl}(A_x), \ l_x = l \}  \ \text{-- which has cardinality at most} \ d_l$. For each $x \in X$, fix an order to the $m_x$ possible ways to sample the $l_x$ variables of vbl($A_x$) generating a bad setting of $\mathscr{Z}_x$. Let $x_{\psi_j}$ be the variable of the smallest index whose value is $\circ$  at the~ $i$-th step of the algorithm. The algorithm then assigns the first entry that has not been utilized from the vector $\Psi_j$ to $x_{\psi_j}$. If doing so makes all variables different from $\circ$ and do not create any bad setting, the algorithm halts. Otherwise, if none bad setting is created, but there are still variables equal to $\circ$, the algorithm goes to the next step. Finally, if some bad setting $(\text{vbl}(A_y),c)$ is created, then fix the $k_y$ variables associated to $\psi_j$ in  $\text{vbl}(A_y)$ and the value $\circ$ is assigned to the  other $l_y$ variables (different from the first $k_y$, note that $\psi_j$ is one of them) of $\text{vbl}(A_y)$. Furthermore, in this case consider a vector $R$ which we assign to its $i$-th entry a tuple $(\alpha, \beta, \gamma)$ where 

\begin{enumerate}

\item $\alpha := l_y$;

\item $\beta \in [d_{l_y}]$ -- that tell us the position of $\text{vbl}(A_y)$ in $\{ \text{vbl}(A_x) : x \in X, \ \psi_j \in \text{vbl}(A_x), \ l_x = l_y \}$;

\item $\gamma \in [m_y]$ -- the label that tell us the values of  the $l_y$ variables of $\text{vbl}(A_y)$.

\end{enumerate}

\medskip

 Then the algorithm goes to the next step. When no bad setting is created, $R_i$ is left empty. The vector $R$ serves as a register since it is essential for the analysis of the algorithm to keep track of what is happening at each step.

\subsection{Algorithm's Analysis}

\

\

$\quad \ \ $Let $X_i$ be the set of variables whose value is $\circ$ after the $i$-th step and $\phi_i$ the configuration of the variables at step $i$ which determines the value of each variable at the $i$-th step. Let $\Psi_1, \ldots, \Psi_n$ be vectors for which the algorithm returns the pair ($R$,$\phi_t$).

\begin{lemma} At each step $i$, the set $X_i$ is uniquely determined by the vector $(R_j)_{j\leq i}$.
\end{lemma}

We omit the proof since the argument is similar to the original approach, see \cite{Aldo, Entropy Compression}.

\begin{lemma} $((R_j)_{j\leq i},\phi_i)$ stores sufficient information to deduce which were the colors used at each step of the algorithm until step $i$. 
\end{lemma}

\textit{Proof.} We use induction. After the first step, $\phi_1$ is the configuration where only the first variable has value different from $\circ$ and is $(\Psi_1)_1$. Suppose that the lemma is true for~ $i - 1$. By the previous lemma we know $X_{i-1}$, in particular, we know the variable~ $x_{\psi_j}$ which received a value at step $i$. In fact, $j= \min \{ s : x_{\psi_s} \in X_{i-1} \}$. We have two situations to analyse: 

Case 1: If $R_i$ is empty, no bad setting was created, then $\phi_{i-1}$ is obtained from $\phi_i$ by simply assigning the value $\circ$ to the variable $x_{\psi_j}$. Hence, by the induction hypothesis, we know which colors were used at each step of the algorithm until step $i-1$. Clearly, the color used at step $i$ is the value assigned to $x_{\psi_j}$ in $\phi_i$. 

Case 2: If $R_i = (\alpha, \beta, \gamma)$ we know which bad setting $(\text{vbl}(A_y),c)$ involving $\psi_j$ was created at step $i$ and which were the values of its~ $l_y$ variables. Thus $\phi_{i-1}$ is obtained from $\phi_i$ by assigning those values to the $l_y$ variables that received the value $\circ$. Then again the colors used at each step until step $i-1$ are determined by the induction hypothesis. By $\gamma$ we know the color assigned to $\psi_j$ at step $i$. \qed

\medskip

From this point, the analysis follows closely the argument in \cite{Entropy Compression} leading to the following variation of the Esperet-Parreau theorem:

\begin{theorem}\label{main} Let $G$ be a graph with vertex set $V$ for which we are seeking a certain vertex-coloring such that there are bad settings to avoid. Let :

\vskip.1cm

1. $\mathcal{Z} := \{ Z \subset V : \text{There exist } c \text{ such that } (Z,c) \text{ is} $ $\text{a bad setting} \} \text{ indexed by a finite}\\ \text{\ set }  I \subset \mathbb{N}$.

\vskip.1cm

2. For each $i \in I$: 

\begin{minipage}[t]{0.063\textwidth}
$\quad \ \ $a.
\end{minipage}\begin{minipage}[t]{\textwidth}
$C_i := \{ c :  (Z_i, c) \ \text{is a bad}$ $\text{setting} \}$. \end{minipage}
\vskip.1cm

\begin{minipage}[t]{0.063\textwidth}
$\quad \ \ $b. \end{minipage}\begin{minipage}[t]{0.83\textwidth}For any vertex $v$ of $Z_i$, let $k_i$ be fixed vertices different from $v$  in $Z_i$ for which, if we know their colors, there are at most $m_i$ ways to extend this partial coloring to a coloring $c \in C_i$. 
\end{minipage}

\vskip.1cm

$\quad \ \ $c. $l_i := \vert Z_i\vert - k_i.$

\vskip.1cm

3. $E := \{ l \in \mathbb{N} : \exists \; i \in I \text{ such that } l_i = l\}$.

\vskip.1cm
\begin{minipage}[t]{0.027\textwidth}
4. \end{minipage}\begin{minipage}[t]{0.866\textwidth}For each $l \in E$ and $v \in V$, let $Q_l(v)$ be the set of subsets of $V$ containing $v$ and equal to some $Z_j \in \mathcal{Z}$ with $l_j = l$. Then $d_l := \max_{v \in V} \vert Q_l(v) \vert$.
\end{minipage}

\vskip.1cm

5. $\phi_E(x) := 1 + \sum_{i\in E}x^i$ and $\gamma := \phi_E^{'}(\tau)$, where $\tau$ satisfies $\phi_E(\tau) - \tau\phi_E^{'}(\tau) = 0$.

\

Then there exists a coloring of $V$ that avoids all bad settings using $\lceil \gamma \sup_{i \in I} (d_{l_i} m_i)^{1/l_i} \rceil$  colors.

\end{theorem}

 \

\textbf{Remark 1.} As explained in \cite{Tao}, to use the entropy compression argument it is essential to compress information losslessly, since in this case there is a bound for how much one can compress the data, for details see \cite{Ash, Shannon}. The difference of this variation in respect to the original approach of Esperet and Parreau \cite{Entropy Compression}, see also \cite{Aldo}, is that we add new information in the register $R$, namely the label $\gamma$ of the partial configuration. The same modification in the register already appears in \cite{Entropy3}.  Since this class of problems in general deals with a finite number of parameters (vertices, colors, etc.), at the beginning we can always fix these labels of the configurations to avoid the second constraint involving the number $k_{(H,c)}$ of the original approach. In this way we can still ensure that the algorithm performs a lossless compression, which was proved in Lemma 2. 

\

\textbf{Remark 2.} Since $d_{l}$ and $m_i$ treat every bad setting equally (giving the same weight to all them), our analysis only works for the uniform case. A more general approach is made in \cite{AIS}, where they recover the result obtained in \cite{GST} using a non-uniform distribution for the $k$-SAT problem.

\section{Legitimate Coloring in Finite Projective Planes}

 Finite projective planes are combinatorial objects with finite elements (lines are finite sets of points) which share properties with classical projective geometry. We denote by $P_n = (P,L)$ projective planes of order $n$, with set of points $P$ and set of lines $L$. For each $n$, $P_n$ has $n^2+n+1$ lines and the same number of points. Every line has $n+1$ points and each point belongs to $n+1$ lines. For precise definitions and properties see the classical reference \cite{Hall}. 

We study the problem of existence of legitimate colorings proposed by Alon and F\"uredi in \cite{Alon}, which can be stated as follows. Let $f$: $P$ $\rightarrow$ $\lbrace 1,\ldots, \mathsf{d} \rbrace$ be a $\mathsf{d}$-coloring of $P$. The \textit{type} of a line $\ell$ with respect to $f$ is defined as the vector $t_{\ell,f} := ( \vert \ell \cap f^{-1}(1) \vert, \ldots , \vert \ell \cap f^{-1}(\mathsf{d}) \vert  )$. If all the lines of a $P_n$ have different types in a $\mathsf{d}$-coloring we say that this is a \textit{legitimate coloring}. When two lines  $\ell_i$ and $\ell_j$ have the same type we say that  $\lbrace \ell_i,\ell_j \rbrace$ is a \textit{bad pair}. The goal is to find the minimum number of colors required - denoted by $\chi(P_n)$ - such that there exists a legitimate coloring of $P_n$ with $\chi(P_n)$ colors. In \cite{Alon}, Alon and  F$\ddot{\text{u}}$redi studied this number for large $n$ using the LLL. More precisely, they proved that if the order of $P_n$ is greater then  $10^{250}$, then  $5 \leqslant \chi(P_n) \leqslant 8$. Furthermore, the authors conjectured that $6 \leqslant \chi(P_n) \leqslant 7$ for large values of $n$, which still is an open problem.

Here we propose a different direction for this question. We investigate what happens when the order of the projective planes are small and more colors are used,  improving the result of Alon and F\"uredi for large $n$ and achieve new results for small values of $n$. In particular, we: prove that $\chi(P_n) \leqslant 8$ when $n \geq 10^{54}$;  guarantee the existence of legitimate colorings for projective planes of any order with 42 colors; and we have results combining the order of planes and the number of colors $\mathsf{d}$. The results are described in Figure 1. 

Our first step is to extend some results obtained in \cite{Alon} from large orders to any order $n$. The results are obtained via a combination of probabilistic method and the entropy compression. 

\\

\begin{center} 
\begin{tikzpicture}[scale = 0.3, every node/.style={scale=0.8}, xscale = 1.8, decoration={
    markings,
   mark=at position 0.5 with {\arrow{>}}}]
    
\draw [help lines, thick, black, ->]  (0,0) -- (23,0);
\draw [help lines, thick, black, ->] (0,0) -- (0,45);
\draw [dashed] (0,0.5) -- (2,0.5);
\node [left] at (0,0.5) {$2$};
\draw [dashed] (0,1.5) -- (3,1.5);
\node [left] at (0,1.5) {$3$};
\draw [dashed] (0,2.5) -- (4,2.5);
\node [left] at (0,2.5) {$4$};
\draw [dashed] (0,3.5) -- (5,3.5);
\node [left] at (0,3.5) {$6$};
\draw [dashed] (0,4.5) -- (6,4.5);
\node [left] at (0,4.5) {$8$};
\draw [dashed] (0,5.5) -- (7,5.5);
\node [left] at (0,5.5) {$11$};
\draw [dashed] (0,6.5) -- (8,6.5);
\node [left] at (0,6.5) {$15$};
\draw [dashed] (0,7.5) -- (9.3,7.5);
\node [left] at (0,7.5) {$20$};
\draw [dashed] (0,9.745) -- (18,9.745);
\node [left] at (0,9.75) {$10^5$};

\draw [dashed] (0,42) -- (2,42);
\node [left] at (0, 44) {$n$};
\node [left] at (0, 42) {$10^{250}$};
\draw [dashed] (0,40.4) -- (3,40.4);
\node [left] at (0, 40.4) {$ 10^{64}$};
\draw [dashed] (0,37.2) -- (2,37.2);
\node [left] at (0, 37.2) {$ 10^{54}$};
\draw [dashed] (0,34) -- (4,34);
\node [left] at (0,34) {$ 10^{37}$};
\draw [dashed] (0,31) -- (5,31);
\node [left] at (0, 31) {$ 10^{27}$};
\draw [dashed] (0,28.2) -- (3,28.2);
\node [left] at (0, 28.2) {$ 10^{25}$};
\draw [dashed] (0,25.6) -- (6,25.6);
\node [left] at (0, 25.6) {$ 10^{22}$};
\draw [dashed] (0,23.2) -- (7,23.2);
\node [left] at (0, 23.2) {$ 10^{19}$};
\draw [dashed] (0,21) -- (8,21);
\node [left] at (0, 21) {$ 10^{17}$};
\draw [dashed] (0,19) -- (9,19);
\node [left] at (0, 19) {$ 10^{15}$};
\draw [dashed] (0,17.2) -- (5,17.2);
\node [left] at (0, 17.2) {$ 10^{13}$};
\draw [dashed] (0,15.6) -- (6,15.6);
\node [left] at (0, 15.6) {$ 10^{11}$};
\node [left] at (0, 14.2) {$ 10^{9}$};
\draw [dashed] (0,14.2) -- (7,14.2);
\node [left] at (0, 13) {$ 10^{8}$};
\draw [dashed] (0,13) -- (8,13);
\node [left] at (0, 12) {$ 10^{7}$};
\draw [dashed] (0,12) -- (9,12);
\node [below] at (2,0) {$8$};
\node [below] at (3,0) {$9$};
\node [below] at (4,0) {$10$};
\node [below] at (5,0) {$11$};
\node [below] at (6,0) {$12$};
\node [below] at (7,0) {$13$};
\node [below] at (8,0) {$14$};
\node [below] at (9,0) {$15$};
\node [below] at (18,0) {$42$};
\node [below] at (22,-0.1) {$\chi(P_n)$};
\path (2,46) coordinate (A) (2,42) coordinate (B) (3,42) coordinate (C)(3,40.4) coordinate (D) (4,40.4) coordinate (E) (4,34) coordinate (F)(5,34) coordinate (G) (5,31) coordinate (H) (6,31) coordinate (I) (6,25.6) coordinate (J) (7,25.6) coordinate (K)(7,23.2) coordinate (L) (8,23.2) coordinate (M) (8,21) coordinate (N) (9,21) coordinate (O) (9,19) coordinate (P) (12.25,19) coordinate (Q) (12.25,17.7) coordinate (R) (15.5,17.7) coordinate (S) (15.5,16.7) coordinate (T) (17.75,16.7) coordinate (U)(17.75,16) coordinate (V) (24, 16) coordinate (X)(24,46) coordinate (Y) ;
  \draw[pattern=north east lines, dashed] (A) -- (B) -- (C) -- (D) -- (E) -- (F) -- (G) -- (H) -- (I) -- (J)-- (K)-- (L)-- (M)-- (N)-- (O)-- (P)-- (Q)-- (R) -- (S) -- (T) -- (U) -- (V) -- (X) -- (Y) -- cycle;
  
  \path  (2,0) coordinate (A)(2,0.5) coordinate (B) (3,0.5) coordinate (C)(3,1.5) coordinate (D) (4,1.5) coordinate (E) (4,2.5) coordinate (F) (5,2.5) coordinate (G) (5,3.5) coordinate (H) (6,3.5) coordinate (I) (6,4.5) coordinate (J) (7,4.5) coordinate (K) (7,5.5) coordinate (L)(8,5.5) coordinate (M)(8,6.5) coordinate (N)(9,6.5) coordinate (O) (9,7.5) coordinate (P) (11.25, 7.5) coordinate (Q) (11.25,8.06) coordinate (R) (13.5,8.06) coordinate (S) (13.5, 8.62) coordinate (T) (15.75,8.62) coordinate (U) (15.75,9.18) coordinate (V) (18,9.18) coordinate (X)(18,9.75) coordinate (Y) (18, 10.31) coordinate (Z) (15.75,10.31) coordinate (a) (15.75,10.87) coordinate (b) (13.5,10.87)  coordinate (c) (13.5,11.43) coordinate (d) (11.25,11.43) coordinate (e)  (11.25,12) coordinate (f) (9,12) coordinate (g) (9,13) coordinate  (h) (8,13) coordinate (i) (8,14.2) coordinate (j) (7,14.2) coordinate (k) (7,15.6) coordinate (l) (6,15.6) coordinate (m) (6,17.2) coordinate (n) (5,17.2) coordinate (o) (5,21) coordinate (p) (4,21) coordinate (q) (4,28.2) coordinate (r) (3,28.2) coordinate (s) (3,37.2) coordinate (t) (2,37.2) coordinate (u)(2,42) coordinate (v) (24,42) coordinate (x) (24,0) coordinate (y);
  \draw[pattern=vertical lines] (A) -- (B) -- (C) -- (D) -- (E) -- (F) -- (G) -- (H) -- (I) -- (J) -- (K) -- (L) -- (M) -- (N) -- (O) -- (P) -- (Q) -- (R) -- (S) -- (T) -- (U) -- (V) -- (X) -- (Y) -- (Z) -- (a) -- (b) -- (c) -- (d) -- (e) -- (f) -- (g) -- (h) -- (i) -- (j) -- (k) -- (l) -- (m) -- (n) -- (o) -- (p) -- (q) -- (r) -- (s) -- (t) -- (u) -- (v) -- (x) -- (y) -- cycle;

\node at (13.5,43.3) {\small\textbf{Previous Region}};
  
\node at (15.3,31) {\small\textbf{LLL Region}};
 \node at (18,5) {\small\textbf{Entropy}};
 \node at (18,4.2) {\small\textbf{Compression}}; 
 \node at (18,3.4) {\small\textbf{Region}};  
 \node[align=center, below, scale=1.2] at (11.5,-2)%
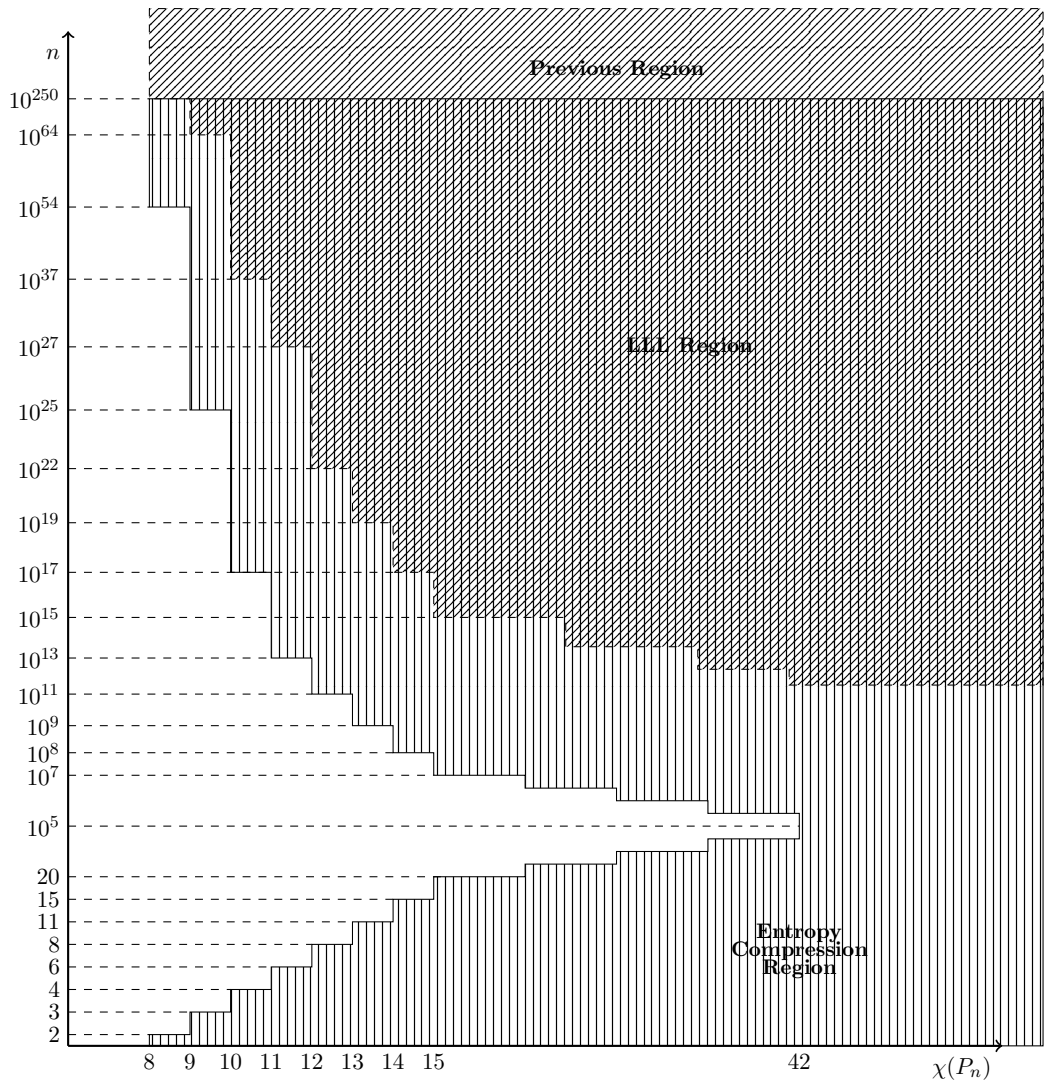
 { Figure 1: New Results for $\chi(P_n)$.};
 \end{tikzpicture}
\end{center}

The next fact is central to the whole approach:

\medskip

\textbf{Fact 1.} \textit{Let $P_n = (P,L)$ be a projective plane of order $n$. Then there exist a set $S \subset P$ such that for all $\ell \in L$  $\frac{1}{2}\ln n   \leqslant   \vert \ell \cap S \vert   \leqslant   11\ln n.$}

 \medskip

\textit{Proof.} In \cite{Alon} is shown that this holds for large $n$ with $20\ln n$ as upper bound instead of $11\ln n.$ We achieved this result applying Corollary 2.8 of \cite{Matas} which is an improvement on Hoeffding-Azuma inequality obtained by McDiarmid in \cite{McDiarmid}.\qed

 \medskip

 From now on $S$ will denote such set.
 
 \newpage
 
 \begin{corollary} 
The cardinality of $S$ is at most $\frac{(n^2+n+1) 11\ln n}{n+1}$.
\end{corollary}

\textit{Proof.} Since there are $n+1$ lines passing through each point of $S$ and $\vert L \vert = n^2 + n + 1$, we have the inequality  $\max_{\ell \in L} \vert \ell \cap S \vert   \geqslant \frac{\vert S \vert (n+1)}{n^2+n+1}$. By Fact 1, for each line $\ell \in L$ we know that $\vert \ell \cap S \vert \leqslant 11\ln n$, then $\vert S \vert \leqslant \frac{(n^2+n+1)11\ln n}{n+1}$.\qed

\medskip

\begin{definition}Let $f \colon P\setminus S \rightarrow \lbrace 1,\ldots, \mathsf{d} \rbrace$ be a partial $\mathsf{d}$-coloring of $P_n$. Given  $\ell_i,\ell_j \in L$, with $\ell_i \neq \ell_j$, we define $\lbrace \ell_i,\ell_j\rbrace $ as a dangerous pair if $d_1(t_{\ell_i,f} , t_{\ell_j,f})$ $\leqslant$ $22\ln n$.
\end{definition}

\medskip

It is not hard to see that if two lines do not form a dangerous pair than there is no extension of the partial coloring $f$ in which these two lines form a bad pair. 

\medskip
 
\textbf{Fact 2.} \textit{Let $\mathsf{d} \geq 8$, $a\geq 1$ and $b\geq 4$. Then, there exists $f \colon P\setminus S \rightarrow \lbrace 1,\ldots, \mathsf{d} \rbrace$ partial $\mathsf{d}$-coloring such that, for $n$ large enough we have:}

 \textit{1. There is no point in $S$ which belongs to more than $b$ lines involved in dangerous pairs.}
 
 \textit{2. There is no line in $L$ which forms a dangerous pair with more than $a$ other lines.}

\medskip

\textit{Proof.} The argument is the same as given in the original paper of Alon and F\"uredi where the authors take $\mathsf{d} = 8, a = 1$ and $b= 4$, see Lemma 2.6, Lemma 2.7 and Proposition 2.8 of \cite{Alon}. The advantage of taking greater values for $a$ and $b$ is that this modification allows us to consider smaller values for $n$, our main contribution to the problem. \qed

\medskip

Now we have all the ingredients to apply the entropy compression method on the problem of legitimate colorings of finite projective planes. Let $P_n = (P,L)$ be a finite projective plane with $P = \{p_1, \ldots, p_{n^2+n+1}\}$ and $L = \{\ell_1, \ldots, \ell_{n^2+n+1}\}$. Fix a set $S \subset P$ and a partial coloring $f \colon P\setminus S \rightarrow [\mathsf{d}]$ satisfying the properties as in the facts 1 and 2. Thus, no point in $S$ belongs to more than $a\cdot b$ dangerous pairs. In order to use Theorem 1, we define our bad settings as follows:

\vskip.1cm

1. $\mathcal{Z} = \{ Z \subset P : \exists  \ \ell_{\kappa},\ell_j, \ell_{\kappa} \neq \ell_j \text{ such that } Z = (\ell_\kappa \cup \ell_j)\backslash (\ell_\kappa \cap \ell_j)  \} \text{ indexed by a finite}\\ \text{\ set }  I \subset \mathbb{N}$.

\vskip.1cm

2. For each $i \in I$: 

\begin{minipage}[t]{0.063\textwidth}
$\quad \ \ $a.
\end{minipage}\begin{minipage}[t]{0.86\textwidth}
 $C_i = \{ c : Z_i = (\ell_\kappa \cup \ell_j)\backslash (\ell_\kappa \cap \ell_j) \ \text{and} \  t_{\ell_\kappa,c} = t_{\ell_j,c} \}.$ \end{minipage}

\vskip.1cm

\begin{minipage}[t]{0.063\textwidth}
$\quad \ \ $b.
\end{minipage}\begin{minipage}[t]{0.83\textwidth}For any point $p$ of $Z_i$, consider the set $Y_i$ containing $p$ and other $\underline{m}-1$ points of smallest index of the line containing $p$, where $\underline{m} = \argmin_m \frac{m}{m-1}(m!a \cdot b (m-1))^{(1/m)}$ (this choice will be clear later). Fixing the points of $Z_i \setminus Y_i$ it is easy to see that $m_i = \underline{m}!$.
\end{minipage}

\vskip.1cm

$\quad \ \ $c. $l_i = \vert Z_i \vert - \vert Z_i \setminus Y_i\vert =  \underline{m}$.

\vskip.1cm

3. $E = \lbrace \underline{m} \rbrace$.

\vskip.1cm

4. $d_{\underline{m}} = a\cdot b$. 
 
\vskip.1cm

\noindent Then, Theorem \ref{main} implies the following result:

\begin{proposition} Let $P_n$ be a projective plane of order $n$. There exists a legitimate coloring of $P_n$ with $\max(\mathsf{d}(n, a, b), \lceil \min_m\frac{m}{m-1}(m!a\cdot b (m-1))^{1/m}\rceil)$ colors, where $\mathsf{d}(n,a,b)$ is the number of colors needed in order to guarantee the existence of a partial coloring $f: P\setminus S \rightarrow [\mathsf{d}(n,a,b)]$  for which no point belongs to more than $a\cdot b$ dangerous pairs.
\end{proposition}

\textit{Proof.} Our first step is to fix a partial coloring  $f$ of $P_n$ for which no point belongs to more than $a\cdot b$ dangerous pairs using $\mathsf{d}(n,a,b)$ colors, we will estimate this number later. Since $E = \lbrace \underline{m} \rbrace$ and $\ \phi_E(x) = 1 + x^{\underline{m}} \ $ we have $ \phi_E(\tau) - \tau\phi_E^{'}(\tau) = 0$ iff  $\tau = (\underline{m} - 1)^{-1 / \underline{m}}$. Then $\phi_E^{'}(\tau) = \underline{m}(\underline{m}-1)^{(1-\underline{m})/\underline{m}}.$
Therefore, using $\lceil \min_m \frac{m}{m-1}(m!a\cdot b (m-1))^{1/m}\rceil$ colors it is possible to extend $f$ to a coloring $C$ such that no bad pairs are created.  \qed

\\

Now, the arguments from \cite{Alon} tell us how to estimate $\mathsf{d}=\mathsf{d}(n,a,b)$, $n$, $a$ and $b$. We do this with more accurate constants. Defining:
$$K := 2^\mathsf{d} \binom{22 \log n + \mathsf{d} }{\mathsf{d}}\frac{\mathsf{d}^{\frac{\mathsf{d}}{2}}}{[2\pi(n + 1 - 11 \log n - n^{1/2})]^{\frac{\mathsf{d}-1}{2}}},$$ a similar argument as in Lemma 2.6 of \cite{Alon} shows that given a line of $L$, the probability $P_a$ of this line to form a dangerous pair with $a + 1$ (instead of 2, as in \cite{Alon})  other lines obeys 
$$ P_a \leq (K)^{a+1} (n^2 + n + 1)\binom{n^2 + n}{a+1}.$$
Moreover, by Lemma 2.7 of \cite{Alon}, given a point $p$ in $P$, the probability $P_b$ of $p$ belonging $b+1$ (instead of 5) lines involved in dangerous pairs with other lines satisfies
$$ P_b \leq \frac{(K)^{b+1}11\log n (n^2 + n + 1)\binom{n + 1}{b+1}(n^2 + n - (b + 1))^{b+1}}{n+1}. $$

\

Thus, the coloring $f$ has the desired properties if $ 1 - (P_a + P_b) > 0$.\\

\medskip

\textbf{Remark 2.} For large values of $n$, the solution of the proposition is given by the following optimization problem:
$$ \text{Minimize} \quad \quad  g(a,b,m,n) = n \quad \quad \quad \quad \quad \quad \quad $$
$$ \text{subject to} \quad \begin{cases}
\frac{m}{m-1}(m!a \cdot b (m-1))^{1/m} \leq \mathsf{d} \\
 1 - (P_a + P_b) > 0

\end{cases} $$

\

which finds the smallest order of a projective plane for which we can find a legitimate coloring using $\mathsf{d}$ colors. Some results are shown in figure 1.

\section{Concluding Remarks}
In this paper, we have shown that every problem formulated in the variable version of the Lov\'asz Local Lemma can be seen as a graph coloring problem where some configurations are forbidden. As a consequence, we show that we can apply the entropy compression method for this class of problems. We presented an application of this situation, the problem of finding legitimate colorings on projective planes. A larger region is obtained where there exists a legitimate coloring with respect to the order of a projective plane including results for small orders. It remains an open problem, though, whether the gap between small and large orders -- see figure 1 -- is intrinsic of the problem or a blind spot of the method used.

\section*{Acknowledgement}

R. Bissacot thanks A. Procacci and Y. Kohayakawa for telling him about the reference \cite{Alon} years ago. Both authors thank G. O. Mota for discussions and references and E. O. Endo, J. Feitoza, R. G. Alves and R. Freire for reading and suggestions in earlier versions of this paper. The authors thank N. Alon for comments about the persistence of the arguments of \cite{Alon} for projective planes of any order and, D. Achlioptas, F. Iliopoulos and A. Sinclair for clarifying the limitation of the algorithm (see Remark 2). L. Doin is supported by FAPESP Grant 15/06444-3 and CNPq Grant 486819/2013-2, R. Bissacot is supported by FAPESP
Grant 2016/25053-8 and CNPq grants 312112/2015-7, 486819/2013-2 and 446658/2014-6.

\end{document}